\documentclass[12pt]{article}
\usepackage{epsfig,amsmath,latexsym,amssymb,color}
\usepackage{graphicx}
\usepackage{amsfonts}
\usepackage{epsfig}

\newtheorem{thm}{Theorem}

\begin{document}
\title{A discretization for the nonlinear parabolic
evolution equation of fractional order in space}
\author{Chien-Hong  Cho\thanks{
Department of Mathematics, National Sun Yat-sen University, 
No.70 Lien-hai Rd., Kaohsiung 804, Taiwan, R.O.C.. Supported by NSTC Grant 113-2115-M-110-006.} ~ and
 Hisashi Okamoto\thanks{
Department of Mathematics, Gakushuin University, Tokyo 
 171-8588 ~ Japan. Partially supported by
JSPS Grant 22K03438  }
}
\maketitle

\begin{abstract}
We consider a nonlinear parabolic equation of fractional order in space 
and propose its numerical discretization.  The fractional derivative  is defined through a functional analytic setting,
rather than the traditional definition of fractional derivatives such as the Riemann-Liouville
derivative.
Numerical experiments are reported and some conjectures are presented. 

\end{abstract}

\section{Introduction}

We consider the following nonlinear 
parabolic equation of fractional order 
\begin{equation}\label{fracPDE0}
\frac{ \partial u}{\partial t} = D^{\gamma} u + f(u),
\end{equation}
where $f(u)$ is a nonlinear term and  $ D^{\gamma}$ is the derivative of order $\gamma$.
Some solutions of this equation exists globally, i.e., in $ 0 < t < \infty$. Some do not 
exist globally and become infinite at finite $t$. If the latter  is the case, we say the solution
blows up. The blow-up problems are popular subjects since in real-world applications blow-up mimics some
 physical phenomena of explosive expansion. 

One of the earliest and the most famous model whose solutions possess finite-time singularities may perhaps be 
the one in Fujita \cite{fujita}. He considered the semilinear heat equation 
with the usual Laplace operator  $\Delta$
\begin{equation}\label{semilinearheat}
\frac{ \partial u}{\partial t}(t,x) = \Delta u(t,x) + u^p(t,x),\quad x\in\mathbb{R}^m,\;t>0, \quad(p>1)
\end{equation}
and showed that the solution may become singular at finite $t$. After his work, 
many researchers devoted themselves to the study of the blow-up solutions of \eqref{semilinearheat}. 
The literatures cannot be listed individually here due to the hugeness of the number of references.
 But one may consult for example \cite{BB,CM,GK,levine,levine2,tsutsumi} for more details. 

In view of the rich results for \eqref{semilinearheat}, researchers have increasingly focused on the following problem 
\begin{equation}\label{fracPDE}
\frac{ \partial u}{\partial t}(t,x) = -(-\Delta)^{\alpha} u(t,x) + u^p(t,x), 
\quad x\in\Omega\subseteq\mathbb{R}^m,\;t>0, 
\end{equation}
where $\alpha\in(0,1)$ and either $\Omega=\mathbb{R}^m$ or $\Omega\subseteq\mathbb{R}^m$ is 
a smooth bounded domain with suitable boundary condition. Similar to the case of $\alpha=1$, 
the solution may become unbounded in finite time. See for instance \cite{FP,Sugitani}.

The purpose of the current paper is to consider the equation \eqref{fracPDE0} and its numerical 
analogue. In Section 2, we review briefly several types of the derivatives of fractional order 
and recall a definition from the viewpoint of functional analysis. Applying this definition, 
we analyze the blow-up behavior of the solutions of the target equation. 
Several numerical experiments are presented in Section 4. From the observation of the 
computational results, some possible behavior is discussed and conjectured in the same section.  

\section{Derivative of fractional order}
In the present paper we consider diffusion of fractional order. Some equations of fractional order are
useful for modeling certain anomalous diffusion: See \cite{BV,OS,OO,OSW}. Also,  J.-L. Lions \cite{lions} 
proposed a modification of the Navier-Stokes equations using a viscosity  of  fractional order.

\subsection{A quick review}
Now, there are many different ways to define $D^{\gamma} $ 
for many different purposes. Perhaps the most well-known definitions
are the Riemann-Liouville operator and the Caputo derivative.
See \cite{D,NOS,TYS} and the references therein. 

For a $2\pi$-periodic function $u(x)$ defined in $ - \pi \le x \le \pi$, it is natural to define
$$
D^{\gamma} u = \sum_{ n \in \mathbb{Z}}  - |n|^{\gamma}  u_n \exp( in x), 
$$
where $ u_n$ is the Fourier coefficients of $ u$:
$$
u(x) = \sum_{n \in \mathbb{Z}} u_n e^{inx}, \qquad  \qquad u_n = \frac{1}{2\pi} \int_{-\pi}^{\pi}  u(x) \exp( -inx). 
$$
Then $ D^2 =  \frac{d^2 u}{dx^2}$. This definition is adopted in the spectral method \cite{STW}. 
Although it is very accurate under some condition, this definition is often costly when we consider 
its numerical discretization.  If we use the Riemann-Liouville operator, there appears another 
problem: it sometimes  does not properly reflect the boundary condition. It uses a particular point,
called a pivot, to define it. For instance, for $ \gamma \in (0,1)$, we define 
$$
D^{\gamma}f(x) = \frac{d}{dx} \frac{1}{\Gamma(\alpha) }  \int_a^x (x-y)^{-\gamma } f(y) dy.
$$
Here $a$ is the pivot. On the other hand, all the points in the interval are equal for the periodic functions
and there is no particular reason to choose $ x = \pm \pi$ or any point in the
interval $ - \pi \le  x \le  \pi$ as a pivot. 
Suppose we have the Dirichlet boundary condition $u=0$ on the boundary. It is not certain whether we
may impose the same boundary condition when $ 0 < \gamma < 1$, particularly when $ \gamma > 0$ is small.

The same comment also applies to the Gr\"unwald-Letnikov derivative, see \cite{OS}. 
It requires us to consider the function in the whole $\mathbb{R}$. 

\subsection{Preliminaries from functional analysis}
We consider it  more natural to define a fractional power of  an operator in the functional analysis 
setting. 
Let $A$ be a self-adjoint positive definite  operator in a Hilbert space $X$. 
Assume for simplicity  that  $A^{-1}$ exists and is a bounded operator in $X$. 
In what follows, we denote by $ {\cal D}(B)$ the domain of definition of the operator
$B$ equipped with the graph norm.

The following definitions of fractional powers  are well-known. 
Suppose that $ 0 < \alpha < 1$ is given. 
Then (see Kato \cite{kato1} for instance) 
\begin{equation}
A^{-\alpha}  = \frac{\sin \pi \alpha}{\pi} \int_0^{\infty} \lambda^{-\alpha} 
(\lambda + A)^{-1} d\lambda.
\label{eq:kato1}
\end{equation}
This is a special case of the following  formula (see\cite{kato1}) for $ \mu \ge 0$:
\begin{equation}
\left( \mu +A^{\alpha} \right)^{-1} = \frac{\sin \pi \alpha}{\pi} \int_0^{\infty}
\frac{ \lambda^{\alpha} }{ \mu^2 + 2 \mu \lambda^{\alpha} \cos \pi \alpha   + \lambda^{2\alpha}  }
(\lambda + A)^{-1} d\lambda.
\label{eq:kato2}
\end{equation}
Indeed, if $ A^{-1}$ exists, we may set $ \mu \downarrow 0$ in \eqref{eq:kato2}
to obtain \eqref{eq:kato1}. 
By the way, \eqref{eq:kato2} is valid if $ \mu > 0$ even in the case where $A^{-1}$ is not bounded.

Now these formulas  imply that
\begin{align}
A^{\alpha} & = 
A^{\alpha - 1} A =  \frac{\sin \pi (1-\alpha)}{\pi} \int_0^{\infty} 
\lambda^{-(1-\alpha)} A   (\lambda + A)^{-1} d\lambda \nonumber \\
& =  \frac{\sin \pi \alpha}{\pi} \int_0^{\infty} 
\lambda^{-(1-\alpha)} \left( 1 - \lambda (\lambda + A)^{-1} \right) d\lambda.
\label{eq:frac2}
\end{align}
This formula \eqref{eq:frac2} is valid for all nonnegative definite operator 
$A$ without assuming the boundedness of $A^{-1}$.
 Its precise meaning is 
that for $ u \in \mathcal{D}(A)$ 
\begin{equation}
A^{\alpha}u 
 =  \frac{\sin \pi\alpha}{\pi} \int_0^{\infty} 
\lambda^{\alpha-1} \left( u - \lambda (\lambda + A)^{-1} u \right) d\lambda.
\label{eq:frac3}
\end{equation}
This is valid also for more general class of operators: i.e., generators of 
analytic semi-groups, see Kato \cite{kato1}.
 We, however, do not use this version in the present paper.

We consider the following fact to be most advantageous: 
\emph{The boundary conditions are included in the domain of the operator} $A^{\alpha}$. 
Accordingly, once $A^{\alpha}$ is defined, we need not be worried about the boundary condition. 
However, the precise determination  of the domain of $ A^{\alpha}$ is rather difficult.
In a special case it can be determined. 
See Fujiwara \cite{fujiwara}. 
Note that we always have $ {\cal D}(A^{\alpha}) \supset {\cal D}(A^{\beta})$ 
for $ 0 < \alpha < \beta \le 1 $. 
Also, it should be noted that our definition works equally even if $A$ is an elliptic differential operator of 
\emph{variable coefficient}.

\medskip

Here is an importnat remark. 
If $ \kappa_0 $ is an eigenvalue of $A$ with its eigenfunction $\phi_0(x)$, 
then $ ( \lambda + A)^{-1} \phi_0(x) = \frac{1}{\lambda + \kappa_0} \phi_0(x)$,
as is verified easily. 
Therefore \eqref{eq:frac2} gives us
$
A^{\alpha} \phi_0(x) = M\phi_0(x),
$
where
\begin{equation}\label{kappa0alpha}
M = \frac{\sin \pi \alpha}{\pi} \int_0^{\infty}
 \lambda^{-(1-\alpha)}  \left( 1 - \frac{\lambda}{ \lambda + \kappa_0} \right) 
 d\lambda = \kappa_0^{\alpha}.
\end{equation}
In fact, one sees that the integral above is equal to 
$
\kappa_0\int_0^\infty\frac{\lambda^{\alpha-1}d\lambda }{\lambda+\kappa_0}=\kappa_0^\alpha\int_0^\infty\frac{t^{\alpha-1}}{1+t}dt
$. 
 Since
$$\int_0^\infty\frac{t^{\alpha-1}}{1+t}dt=B(\alpha,1-\alpha)=\Gamma(\alpha)\Gamma(1-\alpha)=\frac{\pi}{\sin \pi \alpha},$$
we arrive at \eqref{kappa0alpha}. Here $B(\cdot,\cdot)$ denotes the Beta function. Therefore 
\begin{equation}
A^{\alpha} \phi_0 = \kappa_0^{\alpha} \phi_0 \qquad \hbox{if} \qquad
A\phi_0 = \kappa_0 \phi_0.
\label{eq:eigen}
\end{equation}
This can also be proved by what is called the spectral theorem on self-adjoint operators.
But \eqref{eq:eigen} holds true even if $A$ is not self-adjoint but is $m$-accretive.
It is now easy to see that our functional analytic definition produces the same results as the one which we introduced 
via the Fourier series 
in the periodic boundary condition. 

In many applications, the eigenfunction for  the smallest eigenvalue of $A$ is of constant sign throughout
the domain of $x$. Then we can take the same function as an eigenfunction of
$A^{\alpha}$. This fact is useful when we consider applications.

The following remark is rather important. 
If $ (\lambda + A)^{-1}$ preserves positivity, i.e., if $ u \ge 0$ implies
$ ( \lambda + A)^{-1} u \ge 0$ for any $ \lambda > 0$, then $(\lambda + A^{\alpha})^{-1} $ 
also preserves  positivity, as is seen easily by \eqref{eq:kato2}.  Therefore positivity of Green's function of $A$ 
implies positivity of Green's function of $A^{\alpha}$.

Also, we should remark that a functional-analytic approach to evolution equations of fractional order
is presented in \cite{NOS}, which nicely summarizes the results known by 2016. 
Our approach is more specific, and we hope it can be understood more easily. 

\subsection{Examples}
Suppose, for instance, $A$ is $- \frac{d^2}{dx^2}$ in $L^2(\mathbb{R})$.  
Then we obtain  
$$
(\lambda + A) ^{-1} v(x) = \frac{1}{2\sqrt{\lambda}} \int_{-\infty}^{\infty}
e^{- \sqrt{\lambda} |x-y|} v(y) \, dy.
$$
See  \cite{but}, or we can verify its validity by directly differentiating the right hand side. 
If this is plugged into \eqref{eq:frac3}  then we obtain

$$
A^{\alpha}v(x) =  \frac{\sin \pi \alpha}{2\pi} 
\int_{-\infty}^{\infty} \int_0^{\infty}  \lambda^{- \frac{1}{2} + \alpha} 
e^{- \sqrt{\lambda}|x-y|} d\lambda  (v(x) -v(y))dy.
$$
Namely, 
$$
A^{\alpha}v(x) =  \frac{\sin \pi \alpha}{\pi} \Gamma(2\alpha + 1)
\int_{-\infty}^{\infty} \frac{v(x) -v(y)}{|x-y|^{2\alpha + 1}} dy.
$$
The right hand side is well-defined if $ 2\alpha < 1$. If $ 1 \le 2\alpha $, then we proceed
in the following way. 
$$
\int_{-\infty}^{\infty} \frac{v(x) -v(y)}{|x-y|^{2\alpha + 1}} dy=
\int_{-\infty}^{\infty} \frac{v(x) -v(x+z)}{|z|^{2\alpha + 1}} dz=
\int_{-\infty}^{\infty} \frac{v(x) -v(x-z)}{|z|^{2\alpha + 1}} dz.
$$
By addition, we have 
$$
\int_{-\infty}^{\infty} \frac{v(x) -v(y)}{|x-y|^{2\alpha + 1}} dy=
\frac{1}{2}
\int_{-\infty}^{\infty} \frac{2v(x) -v(x+z) - v(x-z)}{|z|^{2\alpha + 1}} dz,
$$
the right hand side of which is well-defined for $ 2\alpha <2$  and a smooth $v$.

In the case of $ A = - \frac{d^2}{dx^2}$, this is not new (see for instance \cite{BV}).
 However, in the case of
general second order operator, our  definition is perhaps more convenient.
For instance it can be applied for defining a fractional power  
of differential equations  of variable coefficients.

In the case of $2\pi$-periodic functions, we have 
$$
(\lambda + A) ^{-1} v(x) = \frac{1}{2\sqrt{\lambda}\mathrm{sinh} \, \pi \sqrt{\lambda}}
 \int_{-\pi}^{\pi}
\mathrm{cosh} \, \left( \sqrt{\lambda} ( \pi - |x-y|) \right)   v(y) \, dy
\qquad ( - \pi < x < \pi).
$$
If the domain is $ 0 \le x \le 1$ and the Dirichlet conditions $u(0)=u(1)=0$ are imposed,
then we have 
$$
(\lambda + A) ^{-1} v(x) = \frac{1}{2\sqrt{\lambda}\mathrm{sinh} \, \sqrt{\lambda}}
 \int_{0}^{1} \left\{  
\mathrm{cosh} \, \left(\sqrt{\lambda} ( 1 - |x-y|) \right)   -
\mathrm{cosh} \, \left(\sqrt{\lambda} ( 1 - x-y) \right)   
\right\}  \, v(y) \, dy
$$
for $ 0 \le x \le 1$. 
If the  Neumenn condition $ u'(0) = u'(1) =0$ are imposed, then  
$$
(\lambda + A) ^{-1} v(x) = \frac{1}{2\sqrt{\lambda}\mathrm{sinh} \, \sqrt{\lambda}}
 \int_{0}^{1} \left\{  
\mathrm{cosh} \, \left(\sqrt{\lambda} ( 1 - |x-y|) \right)   +
\mathrm{cosh} \, \left(\sqrt{\lambda} ( 1 - x-y) \right)   
\right\}  \, v(y) \, dy.
$$
These formulas imply that we can compute $A^{\alpha}$ once we have a means to 
compute  the resolvent $ (\lambda + A)^{-1}$. The resolvent should be discretized by some method
 such as the spectral method or the finite difference method. Once this is done, 
$A^{\alpha}$ is discretized by approximating
the integral \eqref{eq:kato1} or \eqref{eq:frac2}.

\section{Theoretical considerations}

Let us  consider 
\begin{equation}
\frac{\partial u}{\partial t} = -A^{\alpha} u  + u^p, 
\label{eq:ev1}
\end{equation}
where $ p = 2,3,4,\cdots$ and $A$ is a nonnegative definite  self-adjoint operator in 
$L^2 = L^2(\Omega)$. 
Here $\Omega$ denotes some open connected set in $\mathbb{R}^m$.
In order for the solution of \eqref{eq:ev1} to
be constructed, we assume that
\begin{quotation}
(H)~~ $ v \in {\cal D}(A^{\alpha})$ implies 
that $ v^p \in L^2$.
\end{quotation}
 If $A$ is a differential operator of 
second order, Fujiwara's results \cite{fujiwara} 
imply that ${\cal D}(A^{\alpha}) \subset H^{2\alpha}(\Omega)$ where $H^s(\Omega)$ denotes 
the Sobolev space of order $s$. Therefore $ u^p$ has a meaning as an element of $L^2$
for $ u \in {\cal D}(A^{\alpha})$ if $ m = 1$ and $ \alpha > 1/4$. The hypothesis (H) above 
is satisfied if these conditions are assumed. 

\begin{thm}\label{thm1}
If $A$ satisfies the hypothesis \textnormal{(H)}, 
the evolution equation  \eqref{eq:ev1} has a strong solution, local in time, for
all $ u_0 \in {\cal D}(A^{\alpha})$. 
\end{thm}
This theorem can be proved in a standard way, see for instance \cite{BCO}, \cite{OO},  or \cite{pazy}. 
If $ \alpha \le 1/4$, we can  obtain a local solution for 
$ u_0 \in {\cal D}(A^{\beta})$ with suitable $ \beta > \alpha$.  But we do not
pursue the theorem in this direction.

Let $T_{\alpha}$ be the blow-up time of the solutions of \eqref{eq:ev1} with varying $\alpha$ 
with one and the same initial data $u(0)$. It is natural to conjecture that 
$ T_{\alpha} $ is a monotonically decreasing function of $\alpha$ in $(0,1)$, although 
we do not have a proof. See the discussion in Section 4.

We first prove that $T_{\alpha}$  is finite under a certain condition.  
Since we are no longer able to use the maximum principle, we follow the method in 
Levine \cite{levine}. His theorem (Theorem 1 of \cite{levine}) 
reads, in our restricted setting, as follows:
Let $\cal A$ be a nonnegative self-adjoint operator in a Hilbert space $X$. 
The nonlinear term, $f(u) = u^p$, can be more general, but we do not pursue 
a general theorem. We assume that $ u \mapsto u^p$ is a well-defined 
bounded operator from $ {\cal D}({\cal A})$ into $X$. This condition 
is satisfied if we consider the  one-dimensional problem and $
{\cal A} = \left(- \frac{d^2}{dx^2} \right)^{\alpha}$ and $\alpha > 1/4$. 
Let $G$ be defined by 
$
G(u) = \frac{1}{p+1} ( u^p, u),
$
where $ (~,~)$ denotes the inner product of $X$. We now have 
\begin{thm}
If the initial value satisfies 
$$
G\left( u(0) \right) > \frac{1}{2} \left( u(0), {\cal A} u(0) \right), 
$$
then the solution blows up in finite time. 
\end{thm}
The advantage of the use of this theorem is its abstractness. $\cal A$ 
can be any operator, not necessarily a differencial operator, as far as it is self-adjoint and
nonnegative definite. 
We can use this theorem with $ {\cal A} = A^{\alpha}$, where $A$ is the one in \eqref{eq:ev1}.

In applications $A$ is often the minus Laplace operator $- \Delta$ or 
an elliptic operator of second order. In the case of $\alpha=1$, the evolution equation
satisfies the maximum principle. However, if we consider the equation 
with $\alpha\in(0,1)$, the evolution equation does not satisfy
the maximum principle, although positivity is often preserved. 
Therefore many arguments which rely on the maximum principle 
is no longer helpful for the equations of fractional order. 
Let $K$ be a positive cone given by
$$
K = \left\{ f \in L^2(\Omega) \; |  \; f(x) \ge 0 \; \hbox{almost everywhere in}\;  \Omega \right\}.
$$
Suppose that $A$ is positivity preserving in the sense that, for all $ \lambda > 0$, 
$ (\lambda + A )^{-1} K \subset K$. Then  $A^{\alpha}$, too, is positivity preserving. 
This follows immediately from \eqref{eq:kato2}.
Then the semi-group generated by $ -A^{\alpha} $ is
also positivity preserving:
$$
e^{- t A^{\alpha}} K \subset K.
$$
This can be seen from Hille's representation of semi-groups (\cite{kato2}, page 481):
$$
e^{-t{\cal A} } = \lim_{n \rightarrow \infty}  \left( 1 +\frac{t}{n} {\cal A} \right)^{-n}.
$$

Even if $A$ satisfies the maximum principle, $ A^{\alpha}$ does not necessarily  satisfy it. 
Nevertheless positivity preserving property, which is weaker than the maximum principle, 
 is guaranteed. 
Therefore if $ u(x)$ is positive almost everywhere, $ e^{-t A^{\alpha}}u$ 
is positive almost everywhere. Although we cannot use the maximum principle,
this positivity preserving property alone is sometimes of great help.

We now present a second proof of $T_{\alpha} < \infty$ which works in a more restricted situation. 
Here we assume that $ A = - \frac{d^2}{dx^2}$ in 
a bounded domain with the periodic, Dirichlet, or Neumann boundary condition.
In each case we have an everywhere positive eigenfunction  $\phi_0 $ such that 
$ A \phi_0 = \kappa_0 \phi_0$ with $ \kappa_0 \ge 0$. 
Then
$$
\frac{d}{dt} (u(t),\phi_0) = - ( A^{\alpha}u, \phi_0) + (u^p,u)
= - \kappa_0^{\alpha}  (u,\phi_0) + (u^p, \phi_0)
$$
Since $A$ is positivity preserving, as we have seen above, 
we henceforth consider only those $u$'s  which are nonnegative everywhere.
Since $\phi_0$ is so chosen that it is positive everywhere, Jensen's
inequality yields $ (u^p, \phi_0) \ge c (u,\phi_0)^p$, where $c$ is a positive constant. 
Setting $ \varphi(t) = (u,\phi_0)$, we have 
$$
\frac{d}{dt} \varphi(t) \ge - \kappa_0^{\alpha} \varphi(t) + c \varphi(t)^p.
$$
If we assume that $- \kappa_0^{\alpha} \varphi(0) + c \varphi(0)^p > 0$, 
 we can conclude, as is proved in a usual way,   that $\varphi(t)$ becomes infinite in finite time.

\bigskip

Now the following questions naturally arise:

Question 1.  At the blow-up time, is $u$ finite except at a point? This is generically  
true for $ \alpha =1$. But is it also true for $ 0 < \alpha < 1$? This seems to be nontrivial, since
 the operator is no longer a local operator for $ 0 < \alpha < 1$.

Question 2. Is there any numerical method 
which reproduces the property of the blow-up profile?  This may depend on the answer to Question 1. 

Question 3. For a fixed initial data $ u(0)$, is the blow-up time $T_{\alpha}$ 
 a decreasing function of $\alpha$?

Question 4. What modification is necessary for a theorem of Fujita type? 
Will Fujita's exponent change if we replace $-\Delta$ with $ (-\Delta)^{\alpha}$.

\bigskip

Questions 1 and 2 will not be touched in the current paper and Question 3 will be discussed in Section 4.3.
 We here give an answer to Question 4. Consider the semilinear heat equation \eqref{semilinearheat} 
with $ u_0(x) \ge 0$ and $ u_0 \not\equiv 0$. That is, we consider nonnegative solutions of 
\begin{equation}\label{fracPDE:Rm}
\frac{\partial u}{\partial t}(t,x) = -( - \Delta)^{\alpha }u(t,x) + u^p(t,x),\quad x\in\mathbb{R}^m,\;t>0.
\end{equation}
With $\alpha=1$, Fujita \cite{fujita} proved the following 
\begin{thm}
If $  0 < m(p-1) < 2$, then all the solutions blow up in finite time.
If $  2 < m(p-1) $, then all the solutions with sufficiently large initial
data blow up in finite time, solutions having small initial data remain 
smooth for all $t$.
\end{thm}
Apparently the statement is inexact but the reader may grasp the essence of 
the claim. See \cite{fujita} for exact statement. 
The case of the critical exponent, $ m(p-1) = 2$,  was resolved by Hayakawa \cite{hayakawa}
and others:
If  $  2 = m(p-1) $, then the situation is the same as in the case of 
$ m(p-1) < 2$. See \cite{levine2}.

We now consider \eqref{fracPDE:Rm} with $\alpha<1$. Then we can prove 
\begin{thm}
If $  0 < m(p-1) < 2\alpha$, then all the solutions blow up in finite time.
If $  2\alpha  < m(p-1) $, then all the solutions with sufficiently large initial
data blow up in finite time, solutions having small initial data remain 
smooth for all $t$.
\end{thm}
The proof can be attained if we follow almost exactly the argument in Fujita \cite{fujita}. There, 
he used the fundamental solution $H(t,x)$ of $ \partial_t - \Delta$, and its property that
$ 0 < H(t,x) < ct^{-m/2}$, where $c$ is a positive constant independent of $(t,x)$. 
In our case we do not know the exact form of the fundamental solution to
$ \partial_t + (-\Delta)^{\alpha}$.  However, by using the Fourier transform, 
we can prove that the fundamental solution $H_{\alpha}(t,x)$ 
satisfies 
$$
0 <  H_{\alpha}(t,x) < ct^{-m/(2\alpha)}.
$$
And this is sufficient for us to prove this theorem. 
The details are omitted. We also refer the reader to \cite{Sugitani}, in which 
it was proved that the nonnegative solution of \eqref{fracPDE:Rm} blows up if $m(p-1)\leq 2\alpha$.

\section{Numerical results and discussion}
\subsection{Discretization for $t$ variable}
Before going into the details of the  numerical examples, we would like to point out some results
for the discretization of blow-up problems of differentiation of integer order. 
We relied on ideas expounded in  \cite{C1}--\cite{CO3} and \cite{nakagawa}.  
Some other results include \cite{IS,NS,S,TYS}.

If we choose the explicit finite  difference
\begin{equation}
\frac{u(t+\Delta t) - u(t) }{\Delta t} = -A^{\alpha} u(t)  + u(t)^p, 
\label{eq:ev2}
\end{equation}
then what we must do is to compute 
\begin{equation}
u(t+\Delta t) = ( I - \Delta t A^{\alpha} )u(t)  + \Delta t \, u(t)^p
\label{eq:ev2a}
\end{equation}
If we choose the following  implicit finite difference
\begin{equation}
\frac{u(t+\Delta t) - u(t) }{\Delta t} = -A^{\alpha} u(t+\Delta t)  + u(t)^p,
\label{eq:ev3}
\end{equation}
Then 
\begin{equation}
( I + \Delta t A^{\alpha} ) u(t+\Delta t)  = u(t)  + \Delta t  u(t)^p
\quad \hbox{or} \quad 
u(t+\Delta t)  = ( I + \Delta t A^{\alpha} )^{-1} \left[  u(t)  + \Delta t \, u(t)^p 
\right]. 
\label{eq:ev3a}
\end{equation}

To \eqref{eq:ev3a}, we may directly apply \eqref{eq:kato2}. 

\subsection{Numerical experiments}

Consider 
\begin{equation}\label{eq:fracmain}
\frac{\partial u}{\partial t}  = -\left( - \frac{d^2}{dx^2}   \right)^{\alpha } u + u^2,\quad t>0,\;x\in(-\pi,\pi),
\end{equation} 
with the periodic boundary condition. Set  $ A = - \frac{d^2}{dx^2} $ for convenience.  
Let us consider those solutions which are even in $x$.
This is just for the sake of simplicity and the general cases are quite similarly 
treated.  We must compute a matrix version of $ A^{\alpha}$. For this we use the spectral method. 
In general, any discretization can be used. A good survey on computation of fractional powers of matrices
  can be found in \cite{higham}. Recently a very accurate computation was reported by \cite{TSMKZ}.
We, however, use  in this paper a well-established  spectral method leaving those sophisticated methods
above  for future works.

 Let $N$ be a positive integer and let us divide 
$ [-\pi,\pi]$ into $2N$ subintervals, and set  $h = \pi/N$. 
The grid points are:
$$
x_0 = 0, \quad  x_1 = h,  \quad x_2 = 2h, \quad \cdots, \quad x_{N-1} = (N-1)h, \quad x_{N} = \pi
$$
If a function is represented as 
$$
u(x) = \sum_{n=0}^{N} \gamma_n \cos nx.
$$
Then 
$
v(x) = -\left( A^{\alpha} u \right) (x)
$
is computed in the following way. 
If $ n > 0$, then we have 
$$
\gamma_n = \frac{2}{\pi} \int_0^{\pi} u(x) \cos nx dx,
$$
which we approximate as 
$$
\approx  \frac{h}{\pi} 
\left(  u_0 + 2u_1 \cos nh + 2u_2 \cos 2nh + \cdots + 2u_{N-1} \cos (n(N-1)h) + u_N (-1)^n  \right).
$$
If $ n = 0$, we have 
$$
\gamma_0 =
\frac{h}{2\pi} 
\left(  u_0 + 2u_1 + 2u_2  + \cdots + 2u_{N-1} + u_N  \right).
$$
In both cases we applied the trapezoidal rule for the integration. 
We now obtain 
$$
v_m = -(A^{\alpha} u)(x_m) = -\sum_{n=0}^N \gamma_n n^{2\alpha} \cos( n m h)
= -\sum_{n=1}^N \gamma_n n^{2\alpha} \cos( n m h) \qquad (0 \le m \le N).
$$
This is rewritten in the following matrix form: 
$$
v_m = \sum_{k=0}^{N}  T_{m,k}   u_k.
$$
where
$$
T_{m,k} = - \frac{2h}{\pi}  \sum_{n=1}^{N} n^{2\alpha} \cos( mnh) \cos(knh)\qquad \hbox{if}~~~ 0 < k < N.
$$
$$
T_{m,0} = - \frac{h}{\pi}  \sum_{n=1}^{N} n^{2\alpha} \cos( mnh), \qquad 
T_{m,N} = - \frac{h}{\pi}  \sum_{n=1}^{N} n^{2\alpha} \cos( mnh) (-1)^n.
$$
These are messy expressions, but once $N$ is given, they need to be  computed only once before we compute $u$.
Hence $T_{m,k}$ may be regarded as known quantities. 

\subsubsection{Explicit scheme}
The following finite difference (explicit scheme) is employed.
$$
\frac{u_n^{m+1} - u_n^m}{\tau_m } = \sum_{k=0}^N T_{n,k} u_k^m + \left( u_n^m  \right)^2.
$$
$\tau_0 = 0.001$ is fixed and, following \cite{CHO},  we set
$$
\tau_m = \min\left\{  \tau_0, \frac{c}{ \max_n |u_n^m| } \right\}.
$$

\begin{figure*}[htbp]
\begin{center} \leavevmode
\includegraphics[width=8cm]{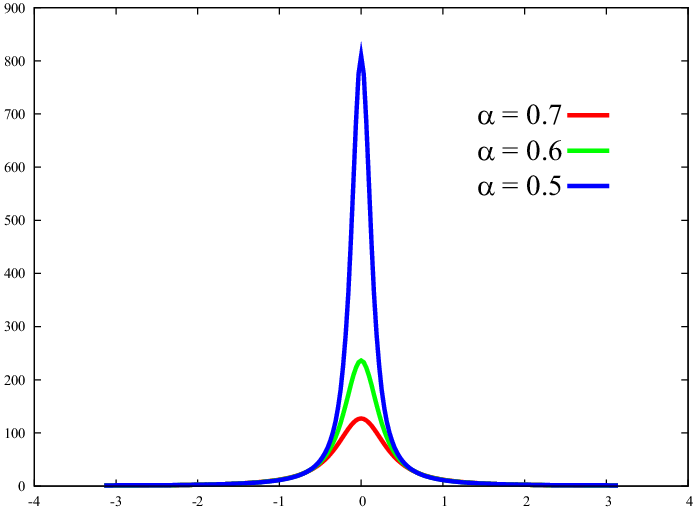}
\end{center}
\caption{ Profiles of $u(0.6,\cdot)$. $\alpha = 0.5,0.6,0.7$ and  $ u(0,x) = \cos x + 1$. 
}\label{zu:time06}
\end{figure*}

We computed solutions for  $ \alpha = 0.5, 0.6,0.7$ with 
$\tau_0 = 0.001, N=100,  u(0,x) = \cos x + 1$. As Figure \ref{zu:time06} shows,  
$\alpha$ certainly works as the smoothing force. And it acts monotonically. 

Note that with these forms, the definition above is valid not only for  
$ 0 < \alpha < 1$ but also for all positive $ \alpha$. However, for a large $\alpha$
we need to choose a very small $\tau_0$. 

Our computation shows that for a fixed $\alpha$ we have $ u(t,x) < u(s,x)$ for $ t < s$. 
We remark that it was proved in Lemma 3.2 in \cite{FP} that the solution of  \eqref{eq:fracmain} 
with zero Dirichlet boundary condition is nondecreasing in $t$ for all $x\in(-\pi,\pi)$ 
if the initial data satisfies 
\begin{equation}\label{initialassump}
-\left(-\frac{d^2}{dx^2}\right)^\alpha u_0+u_0^p\geq 0,\quad\forall x\in(-\pi,\pi).
\end{equation}
The same result also holds for the periodic boundary condition.
 As is easily seen, our initial data $u(0,x)$ satisfies \eqref{initialassump}. In fact, the result can
 be extended, as in the case of $\alpha=1$, as follows: Assume that there exists $t_*\geq 0$ such
 that $$-\left(-\frac{d^2}{dx^2}\right)^\alpha u(t_*,x)+u^p(t_*,x)\geq 0,\quad\forall x\in(-\pi,\pi),$$
then the solution $u$ is nondecreasing in $t$ for all $t\geq t_*$.

\subsubsection{Implicit scheme}
Since $ A^{\alpha}$ is a nonlocal operator, using an explicit scheme is no more
advantageous than implicit schemes. We may therefore consider the following scheme. 
$$
\frac{u(t + \tau_m ) - u(t)}{\tau_m } = - A^{\alpha}  u(t + \tau_m) + \left( u(t)  \right)^2.
$$
i.e., 
$$
\left( 1 + \tau_m A^{\alpha} \right) u(t + \tau_m)  =  u(t) + \tau_m \left( u(t)  \right)^2.
$$

The operator $\left( 1 + \tau A^{\alpha} \right)^{-1} $ may be  defined as follows:
$$
u(x) = \sum_{n=0}^{N} \gamma_n \cos nx, \qquad 
v(x) = \sum_{n=0}^{N} \eta_n \cos nx, \qquad 
\eta_0 = \gamma_0, \qquad  \eta_n = \frac{1}{1 + \tau n^{2\alpha} } \gamma_n
$$
$$
v_m = \left( 1 + \tau A^{\alpha}  \right)^{-1} u(x_m) = 
\sum_{n=0}^N \gamma_n \frac{1}{ 1 + \tau n^{2\alpha}} \cos( n m h)
$$
Accordingly, the matrix representation of $ \left( 1 + \tau A^{\alpha}  \right)^{-1}$ is:
$$
S_{m,k} = \frac{2h}{\pi}  \sum_{n=0}^{N} \frac{\sigma(n)}{1 + \tau n^{2\alpha} }\cos( mnh) \cos(knh)\qquad
 \hbox{if}~~~ 0 < k < N.
$$
$$
S_{m,0} = \frac{h}{\pi}  \sum_{n=0}^{N} \frac{\sigma(n)}{1 + \tau n^{2\alpha} } \cos( mnh), \qquad 
S_{m,N} = \frac{h}{\pi}  \sum_{n=0}^{N} \frac{\sigma(n)}{ 1 + \tau n^{2\alpha} }\cos( mnh) (-1)^n.
$$
where $ \sigma(n)=1 $ for $ 1 < n < N$ and $ \sigma(0)=\sigma(N)= 1/2$.

The advantage of this scheme lies in  its stability. The computational cost is of 
the same order as that of the explicit scheme. 
Indeed, if we use an explicit scheme and $ \alpha $ is large,
 we must take a very small $\Delta t$. That is not the case if
we use the implicit scheme. 
On the other hand, for small $\alpha$, the explicit scheme seems to be more 
stable than the implicit one. 

\begin{figure*}[htbp]
\begin{center} \leavevmode
\includegraphics[width=8cm]{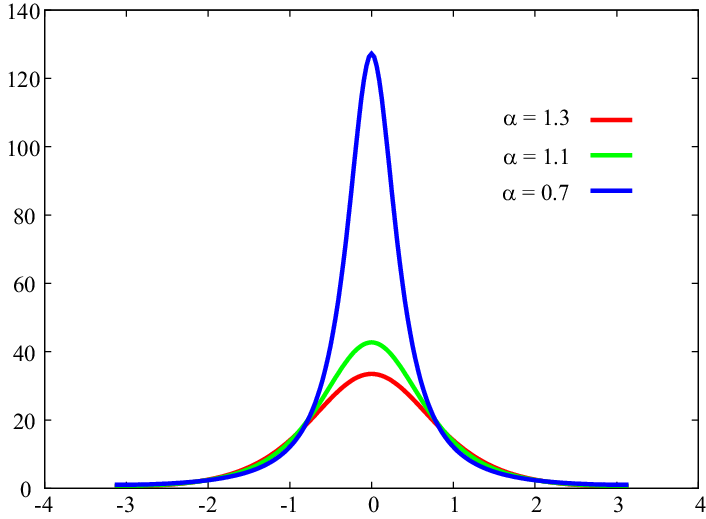}
\end{center}
\caption{   Profiles of $u(0.6,\cdot)$. $\alpha = 0.7,1.1,1.3$ and  $ u(0,x) = \cos x + 1$. 
}\label{zu:time06imp}
\end{figure*}

Moreover, Figure \ref{zu:time06imp}  shows that for a certain $(t,x)$ it is not the case that 
$ u_{\alpha}(t,x) < u_{\beta}(t,x) $ for $ \alpha > \beta$, 
while $ u_{\alpha}(t,0)=\max_x u_\alpha(t,x) < \max_x u_\beta(t,x)=u_{\beta}(t,0) $ holds. 
Here $u_{\alpha}$ and $u_{\beta}$  denote the solution for $ \alpha $ and $ \beta $, respectively.

\subsection{Blow-up time as a function of $\alpha$}
Our computational results suggest that the maximum of the blow-up solution is smaller
 if $\alpha$ is larger. This implies that the blow-up time $T_\alpha$ is nonincreasing 
in  $\alpha$.

Let us consider \eqref{eq:fracmain} with the periodic boundary condition. 
We assume for simplicity that the initial data $ u(0,x)$ is non-negative, even in $x$, and unimodal, whence 
$ u_0(x) := u(0,x)$ satisifes 
\begin{equation}
u_0(x)\ge 0, ~~  u_0'(x) \ge 0, ~~ u_0(x) = u_0(-x) ~~~  for -\pi < x < 0.
\label{assump:IC}
\end{equation}

Then we conjecture that 
\begin{thm}\label{thm:solmono}
Let $0<\beta<\alpha<1$ and $u_\alpha,u_\beta$ be the corresponding solutions
 to \eqref{eq:fracmain}. Suppose that the initial data is not a constant and satisfies
 \eqref{assump:IC}. Then  $u_{\alpha}(t,0) < u_{\beta}(t,0) $, for all $t>0$. 
\end{thm}
\begin{thm}\label{thm:BTmono}
Under the same assumption, 
$T_{\alpha}$, the blow-up time,  is monotonically increasing in $ \alpha$. 
\end{thm}

Theorem \ref{thm:BTmono} is a direct consequence of Theorem \ref{thm:solmono}.
 We encountered some difficulty in proving Theorem \ref{thm:solmono}, but we consider it is very likely to be true.

\end{document}